\documentclass{elsart}
\usepackage{graphics}
\usepackage{graphicx}
\usepackage{epsfig}
\usepackage{amssymb,amsmath}
\newtheorem{Lemma}{Lemma}[section]

\begin{document}
\begin{frontmatter}
\title{Lyapunov Computational Method for Two-Dimensional Boussinesq Equation}
\author{Anouar Ben Mabrouk et al}
\ead{anouar.benmabrouk@issatso.rnu.tn}
\address{Computational Mathematics Laboratory, Department of Mathematics,
Faculty of Sciences, 5000 Monastir, Tunisia.}
\begin{abstract}
A numerical method is developed leading to Lyapunov operators to
approximate the solution of two-dimensional Boussinesq equation.
It consists of an order reduction method and a finite difference
discretization. It is proved to be uniquely solvable and analyzed
for local truncation error for consistency. The stability is
checked by using Lyapunov criterion and the convergence is ......
Some numerical implementations are provided at the end of the
paper to validate the theoretical results.
\end{abstract}
\begin{keyword}
Bousinesq equation, Finite-difference scheme, Stability analysis,
Consistency, Convergence, Error estimates, Lyapunov
operator.\\
\PACS 35B05, 65M06.
\end{keyword}
\end{frontmatter}
\section{Introduction}
The present paper is devoted to the development of a numerical
method based on two-dimensional finite difference scheme to
approximate the solution of the nonlinear Boussinesq equation in
$\mathbb{R}^2$ written on the form
\begin{equation}\label{eqn1-1}
u_{tt}=\Delta\,u+u_{xxxx}+(u^2)_{xx},\quad((x,y),t)\in\Omega\times(t_0,+\infty)
\end{equation}
with initial conditions
\begin{equation}\label{eqn1-2}
u(x,y,t_0)=u_0(x,y)\quad\hbox{and}\quad\displaystyle\frac{\partial\,u}{\partial\,t}(x,y,t_0)=\varphi(x,y),\quad(x,y)\in\Omega
\end{equation}
and boundary conditions
\begin{equation}\label{eqn1-3}
\displaystyle\frac{\partial\,u}{\partial\,n}(x,y,t)=0,\quad
((x,y),t)\in\partial\Omega\times(t_0,+\infty)
\end{equation}
on a rectangular domain $\Omega=]L_0,L_1[\times]L_0,L_1[$ in
$\mathbb{R}^2$. $t_0\geq0$ is a real parameter fixed as the
initial time, $u_{tt}$ is the second order partial derivative in
time,
$\Delta=\frac{\partial^2}{\partial\,x^2}+\frac{\partial^2}{\partial\,y^2}$
is the Laplace operator in $\mathbb{R}^2$, $u_{xx}$ and $u_{xxxx}$
are respectively the second order and the fourth order partial
derivative according to $x$.
$\displaystyle\frac{\partial}{\partial n}$ is the outward normal
derivative operator along the boundary $\partial\Omega$. Finally,
$u$, $u_0$ and $\varphi$ are real valued functions with $u_0$ and
$\varphi$ are $\mathcal{C}^2$ on $\overline{\Omega}$.

The Boussinesq equation
.............................................................\\
................................................................

Our main idea consists in developing a numerical scheme to
approximate the solution of (\ref{eqn1-1})-(\ref{eqn1-3}). The
time and the space partial derivatives are replaced by
finite-difference approximations in order to transform the initial
boundary-value problem (\ref{eqn1-1})-(\ref{eqn1-3}) into a linear
algebraic system. An order reduction method is adapted leading to
a system of coupled PDEs which is transformed by the next to a
discrete algebraic system. The resulting method is analyzed for
local truncation error and stability and the scheme is proved to
be uniquely solvable and convergent.
\section{Two-Dimensional Finite Difference Scheme}
Consider the rectangular domain
$\Omega=]L_0,L_1[\times]L_0,L_1[\subset\mathbb{R}^2$ and an
integer $J\in\mathbb{N}^*$. Denote
$h=\displaystyle\frac{L_1-L_0}{J}$ for the space step,
$x_j=L_0+jh$ and $y_m=L_0+mh$ for all
$(j,m)\in{I}^2=\{0,1,\dots,J\}^2$. Let $l=\Delta\,t$ be the time
step and $t_n=t_0+nl$, $n\in\mathbb{N}$ for the discrete time
grid. For $(j,m)\in{I}$ and $n\geq0$, $u_{j,m,k}^n$ will be the
net function $u(x_j,y_m,t_n)$ and $U_{j,m}^n$ the numerical
solution. The following discrete approximations will be applied
for the different differential operators involved in the problem.
For time derivatives, we set
$$
\displaystyle\,u_t\rightsquigarrow\displaystyle\frac{U_{j,m}^{n+1}-U_{j,m}^{n-1}}{2l}\quad\hbox{and}\quad
\displaystyle\,u_{tt}\rightsquigarrow\displaystyle\frac{U_{j,m}^{n+1}-2U_{j,m}^{n}+U_{j,m}^{n-1}}{l^2}
$$
and for space derivatives, we shall use
$$
\displaystyle\,u_x\rightsquigarrow\displaystyle\frac{U_{j+1,m}^{n}-U_{j-1,m}^{n}}{2h}\quad\hbox{and}\quad
\displaystyle\,u_y\rightsquigarrow\displaystyle\frac{U_{j,m+1}^{n}-U_{j,m-1}^{n}}{2h}
$$
for first order derivatives and
$$
\displaystyle\,u_{xx}\rightsquigarrow\displaystyle\frac{U_{j+1,m}^{n,\alpha}-2U_{j,m}^{n,\alpha}+U_{j-1,m}^{n,\alpha}}{h^2},\quad
\displaystyle\,u_{yy}\rightsquigarrow\displaystyle\frac{U_{j,m+1}^{n,\alpha}-2U_{j,m}^{n,\alpha}+U_{j,m-1}^{n,\alpha}}{h^2}
$$
for second order ones, where for $n\in\mathbb{N}^*$ and
$\alpha\in\mathbb{R}$,
$$
u^{n,\alpha}=\alpha\,U^{n+1}+(1-2\alpha)U^{n}+\alpha\,U^{n-1}.
$$
Finally, we denote $\sigma=\displaystyle\frac{l^2}{h^2}$ and
$\delta=\displaystyle\frac{1}{h^2}$.
\section{Discrete Two-Dimensional Boussinesq Equation}
The method is based on an order reduction idea in which we set
\begin{equation}\label{eqn2-1}
v=u_{xx}+u^2.
\end{equation}
For $(j,m)\in\mathring{I}^2$ an interior point of the grid $I^2$,
($\mathring{I}=\{1,2,\dots,J-1\}$), and $n\geq1$, the following
discrete equation is deduced from (\ref{eqn1-1}).
\begin{equation}\label{eqn1-1discrete}
\begin{matrix}&&U_{j,m}^{n+1}-2U_{j,m}^{n}+U_{j,m}^{n-1}\hfill\cr
&=&\sigma\alpha\left(U_{j-1,m}^{n+1}-4U_{j,m}^{n+1}+U_{j+1,m}^{n+1}
+U_{j,m-1}^{n+1}+U_{j,m+1}^{n+1}\right)\hfill\cr
&&+\sigma(1-2\alpha)\left(U_{j-1,m}^{n}-4U_{j,m}^{n}+U_{j+1,m}^{n}
+U_{j,m-1}^{n}+U_{j,m+1}^{n}\right)\hfill\cr
&&+\sigma\alpha\left(U_{j-1,m}^{n-1}-4U_{j,m}^{n-1}+U_{j+1,m}^{n-1}
+U_{j,m-1}^{n-1}+U_{j,m+1}^{n-1}\right)\hfill\cr
&&+\sigma\alpha\left(V_{j-1,m}^{n+1}-2V_{j,m}^{n+1}+V_{j+1,m}^{n+1}\right)\hfill\cr
&&+\sigma(1-2\alpha)\left(V_{j-1,m}^{n}-2V_{j,m}^{n}+V_{j+1,m}^{n}\right)\hfill\cr
&&+\sigma\alpha\left(V_{j-1,m}^{n-1}-2V_{j,m}^{n-1}+V_{j+1,m}^{n-1}\right).\hfill
\end{matrix}
\end{equation}
Similarly, the following discrete equation is obtained from
equation (\ref{eqn2-1}).
\begin{equation}\label{eqn2-1discrete}
\begin{matrix}\,V_{j,m}^{n+1}+V_{j,m}^{n-1}
&=&2\delta\alpha\left(U_{j-1,m}^{n+1}-2U_{j,m}^{n+1}+U_{j+1,m}^{n+1}\right)\hfill\cr
&&+2\delta(1-2\alpha)\left(U_{j-1,m}^{n}-2U_{j,m}^{n}+U_{j+1,m}^{n}\right)\hfill\cr
&&+2\delta\alpha\left(U_{j-1,m}^{n-1}-2U_{j,m}^{n-1}+U_{j+1,m}^{n-1}\right)\hfill\cr
&&+2\widehat{F(U_{j,m}^{n})}\hfill
\end{matrix}
\end{equation}
where $F(u)=u^2$, $F^n=F(u^n)$ and
$\widehat{F^n}=\displaystyle\frac{F^{n-1}+F^{n}}{2}$. The discrete
boundary conditions are written for $n\geq0$ as
\begin{equation}\label{CBD1}
\displaystyle\,U_{1,m}^{n}=U_{-1,m}^{n}\quad\hbox{and}\quad
\displaystyle\,U_{J-1,m}^{n}=U_{J+1,m}^{n},
\end{equation}
\begin{equation}\label{CBD2}
\displaystyle\,U_{j,1}^{n}=U_{j,-1}^{n}\quad\hbox{and}\quad
\displaystyle\,U_{j,J-1}^{n}=U_{j,J+1}^{n},
\end{equation}
Next, as it is motioned in the introduction, the idea consists in
applying Lyapunov operators to approximate the solution of the
continuous problem (\ref{eqn1-1})-(\ref{eqn1-3}) or its discrete
equivalent system (\ref{eqn1-1discrete})-(\ref{CBD2}). Denote
$$
a_1=\displaystyle\frac{1}{2}+2\alpha\sigma,\quad\,a_2=-\alpha\sigma,
$$
$$
b_1=1-2(1-2\alpha)\sigma,\quad\,b_2=(1-2\alpha)\sigma,
$$
$$
c_1=(1-2\alpha)\delta\quad\hbox{and}\quad\,c_2=\alpha\delta.
$$
Equation (\ref{eqn1-1discrete}) becomes
\begin{equation}\label{eqn1-1discrete-forme-a}
\begin{matrix}&&a_2U_{j-1,m}^{n+1}+a_1U_{j,m}^{n+1}+a_2U_{j+1,m}^{n+1}
+a_2U_{j,m-1}^{n+1}+a_1U_{j,m}^{n+1}+a_2U_{j,m+1}^{n+1}\hfill\cr
&&+a_2\left(V_{j-1,m}^{n+1}-2V_{j,m}^{n+1}+V_{j+1,m}^{n+1}\right)\hfill\cr
&=&b_2U_{j-1,m}^{n}+b_1U_{j,m}^{n}+b_2U_{j+1,m}^{n}+b_2U_{j,m-1}^{n}+b_1U_{j,m}^{n}+b_2U_{j,m+1}^{n}\hfill\cr
&&-a_2U_{j-1,m}^{n-1}-a_1U_{j,m}^{n-1}-a_2U_{j+1,m}^{n-1}-a_2U_{j,m-1}^{n-1}-a_1U_{j,m}^{n-1}-a_2U_{j,m+1}^{n-1}\hfill\cr
&&+b_2\left(V_{j-1,m}^{n}-2V_{j,m}^{n}+V_{j+1,m}^{n}\right)\hfill\cr
&&-a_2\left(V_{j-1,m}^{n-1}-2V_{j,m}^{n-1}+V_{j+1,m}^{n-1}\right).\hfill
\end{matrix}
\end{equation}
Equation (\ref{eqn2-1discrete}) becomes
\begin{equation}\label{eqn2-1discrete-forme-a}
\begin{matrix}
&&V_{j,m}^{n+1}-2c_2\left(U_{j-1,m}^{n+1}-2U_{j,m}^{n+1}+U_{j+1,m}^{n+1}\right)\hfill\cr
&=&2c_1\left(U_{j-1,m}^{n}-2U_{j,m}^{n}+U_{j+1,m}^{n}\right)\hfill\cr
&&+2c_2\left(U_{j-1,m}^{n-1}-2U_{j,m}^{n-1}+U_{j+1,m}^{n-1}\right)\hfill\cr
&&-V_{j,m}^{n-1}+2\widehat{F(U_{j,m}^{n})}.\hfill
\end{matrix}
\end{equation}
Denote $A$, $B$ and $R$ the matrices defined by
$$
A=\begin{pmatrix}a_1&2a_2&0&...&...&0\hfill\cr
a_2&a_1&a_2&\ddots&\ddots&\vdots\hfill\cr
0&\ddots&\ddots&\ddots&\ddots&\vdots\hfill\cr
\vdots&\ddots&\ddots&\ddots&\ddots&0\hfill\cr
\vdots&\ddots&\ddots&a_2&a_1&a_2\hfill\cr
0&...&...&0&2a_2&a_1\hfill
\end{pmatrix},\quad
B=\begin{pmatrix}b_1&2b_2&0&...&...&0\hfill\cr
b_2&b_1&b_2&\ddots&\ddots&\vdots\hfill\cr
0&\ddots&\ddots&\ddots&\ddots&\vdots\hfill\cr
\vdots&\ddots&\ddots&\ddots&\ddots&0\hfill\cr
\vdots&\ddots&\ddots&b_2&b_1&b_2\hfill\cr
0&...&....&0&2b_2&b_1\hfill
\end{pmatrix}
$$
and
$$
R=\begin{pmatrix}-2&2&0&...&...&0\hfill\cr
1&-2&1&\ddots&\ddots&\vdots\hfill\cr
0&\ddots&\ddots&\ddots&\ddots&\vdots\hfill\cr
\vdots&\ddots&\ddots&\ddots&\ddots&0\hfill\cr
\vdots&\ddots&\ddots&1&-2&1\hfill\cr 0&...&...&0&2&-2\hfill
\end{pmatrix}.
$$
The system (\ref{CBD1})-(\ref{eqn2-1discrete-forme-a}) can be
written on the matrix form
\begin{equation}\label{matrixform1}
\left\{\begin{matrix}
\mathcal{L}_A(U^{n+1})+a_2RV^{n+1}=\mathcal{L}_B(U^{n})-\mathcal{L}_A(U^{n-1})+R(b_2V^{n}-a_2V^{n-1}),\hfill\cr
\displaystyle\frac{1}{2}V^{n+1}-c_2RU^{n+1}\!=\!R(c_1U^{n}+c_2U^{n-1})
-\displaystyle\frac{1}{2}V^{n-1}\!+\!\widehat{F^n}\hfill\end{matrix}\right.
\end{equation}
for all $n\geq1$ where
$$
U^n=(U_{j,m}^n)_{0\leq\,j,m\leq\,J},\;V^n=(V_{j,m}^n)_{0\leq\,j,m\leq\,J}\;\hbox{and}\;
F^n=(F(U_{j,m}^n))_{0\leq\,j,m\leq\,J}
$$
and for a matrix $Q\in\mathcal{M}_{(J+1)^2}(\mathbb{R})$,
$\mathcal{L}_Q$ is the Lyapunov operator defined by
$$
\mathcal{L}_Q(X)=QX+XQ,\;\forall\,X\in\mathcal{M}_{(J+1)^2}(\mathbb{R}).
$$
\section{Solvability of the Discrete Problem}
We are now ready to set our first main result.
\begin{thm}\label{theorem1}
The system (\ref{matrixform1}) is uniquely solvable whenever $U^0$
and $U^1$ are known.
\end{thm}
{\it Proof.} It reposes on the inverse of Lyapunov operators.
Consider the endomorphism $\Phi$ which associates to any element
$(X,Y)$ in
$\mathcal{M}_{(J+1)^2}(\mathbb{R})\times\mathcal{M}_{(J+1)^2}(\mathbb{R})$
its image
$\Phi(X,Y)=(AX+XA+a_2RY,\displaystyle\frac{1}{2}Y-c_2RX)$. To
prove Theorem \ref{theorem1}, it suffices to show that $ker\Phi$
is reduced to $0$. Indeed,
$$
\Phi(X,Y)=0\Longleftrightarrow(AX+XA+a_2RY,\displaystyle\frac{1}{2}Y-c_2RX)=(0,0)
$$
or equivalently,
$$
Y=2c_2RX\quad\hbox{and}\quad\,(A+2a_2c_2R^2)X+XA=0.
$$
So, the problem is transformed to the resolution of a Lyapunov
type equation of the form
\begin{equation}\label{Lyapunovequation}
\mathcal{L}_{W,A}(X)=WX+XA=0
\end{equation}
where $W$ is the matrix given by $W=A+2a_2c_2R^2$. Denoting
$$
\omega=2a_2c_2,\,\,\omega_1=a_1+6\omega,\;\;\overline{\omega}_1=\omega_1+\omega
\quad\hbox{and}\quad\omega_2=a_2-4\omega
$$
the matrix $W$ is explicitly given by
$$
W=\begin{pmatrix}\omega_1&2\omega_2&2\omega&0&\dots&\dots&\dots&0&\hfill\cr
\omega_2&\overline{\omega}_1&\omega_2&\omega&\ddots&\ddots&\ddots&\vdots&\hfill\cr
\omega&\omega_2&\omega_1&\omega_2&\omega&\ddots&\ddots&\vdots&\hfill\cr
0&\ddots&\ddots&\ddots&\ddots&\ddots&\ddots&\vdots&\hfill\cr
\vdots&\ddots&\ddots&\ddots&\ddots&\ddots&\ddots&0&\hfill\cr
\vdots&\ddots&\ddots&\omega&\omega_2&\omega_1&\omega_2&\omega&\hfill\cr
\vdots&\ddots&\ddots&\ddots&\omega&\omega_2&\overline{\omega}_1&\omega_2&\hfill\cr
0&\dots&\dots&\dots&0&2\omega&2\omega_2&\omega_1&\hfill
\end{pmatrix}.
$$
Next, we use the following result.
\begin{Lemma}\label{LemmeInversion}
Let $E$ be a finite dimensional ($\mathbb{R}$ or $\mathbb{C}$)
vector space and $(\Phi_n)_n$ be a sequence of endomorphisms
converging uniformly to an invertible endomorphism $\Phi$. Then,
there exists $n_0$ such that, for any $n\geq\,n_0$, the
endomorphism $\Phi_n$ is invertible.
\end{Lemma}
Assume now that $l=o(h^{2+s})$, with $s>0$ which is always
possible. Then, the coefficients appearing in $A$ and $W$ will
satisfy as $h\longrightarrow0$ the following.
$$
A_{i,i}=\displaystyle\frac{1}{2}+\varepsilon\,h^{2+2s}\longrightarrow\displaystyle\frac{1}{2}.
$$
For $1\leq\,i\leq\,J-1$,
$$
A_{i,i-1}=A_{i,i+1}=\displaystyle\frac{A_{0,1}}{2}=\displaystyle\frac{A_{J,J-1}}{2}=-\varepsilon\,h^{2+2s}\longrightarrow0.
$$
For $2\leq\,i\leq\,J-2$,
$$
W_{i,i}=W_{0,0}=W_{J,J}=\displaystyle\frac{1}{2}+2\alpha\varepsilon\,h^{2+2s}-12\alpha^2\varepsilon\,h^{2s}\longrightarrow\displaystyle\frac{1}{2}.
$$
Similarly,
$$
W_{1,1}=W_{J-1,J-1}=\displaystyle\frac{1}{2}+2\alpha\varepsilon\,h^{2+2s}-14\alpha^2\varepsilon\,h^{2s}\longrightarrow\displaystyle\frac{1}{2}
$$
and
$$
W_{i,i-1}=W_{i,i+1}=\displaystyle\frac{W_{0,1}}{2}=\displaystyle\frac{W_{J,J-1}}{2}=-\alpha\varepsilon\,h^{2+2s}+8\alpha^2\varepsilon\,h^{2s}\longrightarrow0
$$
Finally,
$$
W_{i,i-2}=W_{i,i+2}=\displaystyle\frac{W_{0,2}}{2}=\displaystyle\frac{W_{J,J-2}}{2}=-2\alpha^2\varepsilon\,h^{2s}\longrightarrow0.
$$
Next, observing that for all $X$ in the space
$\mathcal{M}_{(J+1)^2}(\mathbb{R})\times\mathcal{M}_{(J+1)^2}(\mathbb{R})$,
$$
\begin{matrix}
\|(\mathcal{L}_{W,A}-I)(X)\|&=\|(W-\frac{1}{2}I)X+X(A-\frac{1}{2}I)\|\hfill\cr
&\leq\left[\|W-\frac{1}{2}I\|+\|A-\frac{1}{2}I\|\right]\|X\|,\hfill
\end{matrix}
$$
it results that
\begin{equation}\label{LWAtendsVersId}
\|\mathcal{L}_{W,A}-I\|\leq\|W-\displaystyle\frac{1}{2}I\|+\|A-\displaystyle\frac{1}{2}I\|\leq\,C(\alpha)h^{2s}.
\end{equation}
Consequently, the Lyapunov endomorphism $\mathcal{L}_{W,A}$
converges uniformly to the identity $I$ as $h$ goes towards 0 and
$l=o(h^{2+s})$ with $s>0$. Using Lemma \ref{LemmeInversion}, the
operator $\mathcal{L}_{W,A}$ is invertible for $h$ small enough.
\section{Convergence of the Discrete Method}
..............

\section{Consistency and Stability of the Discrete Method}
The consistency of the proposed method is done by evaluating the
the local truncation error arising from the discretization of the
system
\begin{equation}\label{system1}
\left\{\begin{matrix}\,u_{tt}-\Delta\,u-v_{xx}=0,\hfill\cr
v=qu_{xx}+u^2.\hfill\end{matrix}\right.
\end{equation}
The principal part of the first equation is
$$
\begin{matrix}\mathcal{L}_{u,v}^1(t,x,y)&=&\displaystyle\frac{l^2}{12}\displaystyle\frac{\partial^4u}{\partial\,t^4}
-\displaystyle\frac{h^2}{12}\left(\displaystyle\frac{\partial^4u}{\partial\,x^4}+\displaystyle\frac{\partial^4u}{\partial\,y^4}\right)
-\alpha\,l^2\displaystyle\frac{\partial^2(\Delta\,u)}{\partial\,t^2}\hfill\cr
&-&\displaystyle\frac{h^2}{12}\displaystyle\frac{\partial^2v}{\partial\,x^4}
-\alpha\,l^2\displaystyle\frac{\partial^4v}{\partial\,t^2\partial\,x^2}+O(l^2+h^2).\hfill\end{matrix}
$$
The principal part of the local error truncation due to the second
part is
$$
\begin{matrix}\mathcal{L}_{u,v}^2(t,x,y)&=&\displaystyle\frac{l^2}{2}\displaystyle\frac{\partial^2v}{\partial\,t^2}
+\displaystyle\frac{l^4}{24}\displaystyle\frac{\partial^4v}{\partial\,t^4}
-\displaystyle\frac{h^2}{12}\displaystyle\frac{\partial^4u}{\partial\,x^4}\hfill\cr
&-&\alpha\,l^2\displaystyle\frac{\partial^4u}{\partial\,t^2\partial\,x^2}+O(l^2+h^2).\hfill\end{matrix}
$$
It is clear that the two operators $\mathcal{L}_{u,v}^1$ and
$\mathcal{L}_{u,v}^2$ tend toward 0 as $l$ and $h$ tend to 0,
which ensures the consistency of the method. Furthermore, the
method is consistent with an order 2 in time and space.

We now proceed by proving the stability of the method by applying
the Lyapunov criterion. A linear system
$\mathcal{L}(x_{n+1},x_{n},x_{n-1},\dots)=0$ is stable in the
sense of Lyapunov if for any bounded initial solution $x_{0}$ the
solution $x_{n}$ remains bounded for all $n\geq0$. Here, we will
precisely prove the following result.
\begin{lem}\label{LyapunovStabilityLemma}
$\mathcal{P}_n$: The solution $(U^n,V^n)$ is bounded independently
of $n$ whenever the initial solution $(U^0,V^0)$ is bounded.
\end{lem}
We will proceed by recurrence on $n$. Assume firstly that
$\|(U^0,V^0)\|\leq\eta$ for some $\eta$ positive. Using the system
(\ref{matrixform1}), we obtain
\begin{equation}\label{LyapunovStability2}
\left\{\begin{matrix}
\mathcal{L}_{W,A}(U^{n+1})=\mathcal{L}_{\widetilde{B},B}(U^{n})+b_2RV^n-\mathcal{L}_{W,A}(U^{n-1})-a_2R(F^{n-1}+F^n),\hfill\cr
V^{n+1}=2c_2RU^{n+1}+2R(c_1U^{n}+c_2U^{n-1})-V^{n-1}+2\widehat{F^n}.\hfill\end{matrix}\right.
\end{equation}
where $\widetilde{B}=B-2a_2c_1R^2$. Consequently,
\begin{equation}\label{LyapunovStability3}
\begin{matrix}
\|\mathcal{L}_{W,A}(U^{n+1})\|\leq\|\mathcal{L}_{\widetilde{B},B}\|.\|U^{n}\|+2|b_2|.\|V^n\|\hfill\cr
\qquad\qquad\qquad\qquad\qquad
+\|\mathcal{L}_{W,A}\|.\|U^{n-1}\|+2|a_2|(\|F^{n-1}\|+\|F^n\|)\hfill\end{matrix}
\end{equation}
and
\begin{equation}\label{LyapunovStability4}
\begin{matrix}
\|V^{n+1}\|\leq4|c_2|.\|U^{n+1}\|+4(|c_1|.\|U^{n}\|+|c_2|.\|U^{n-1}\|)\hfill\cr
\qquad\qquad\qquad\qquad\qquad
+\|V^{n-1}\|+\|F^{n-1}\|+\|F^n\|.\hfill\end{matrix}
\end{equation}
Next, recall that, for $l=o(h^{s+2})$ small enough, $s>0$, we have
$$
a_1=\displaystyle\frac{1}{2}+2\alpha\,h^{2s+2}\rightarrow\displaystyle\frac{1}{2},\quad\,a_2=-\alpha\,h^{2s+2}\rightarrow0,
$$
$$
b_1=1-2(1-2\alpha)h^{2s+2}\rightarrow1,\quad\,b_2=(1-2\alpha)h^{2s+2}\rightarrow0,
$$
$$
c_1=(1-2\alpha)h^{-2}\rightarrow\infty\quad\hbox{and}\quad\,c_2=\alpha\,h^{2s+2}\rightarrow\infty,
$$
$$
a_2c_1=-\alpha(1-2\alpha)h^{2s}\rightarrow0.
$$
As a consequence, for $h$ small enough,
\begin{equation}\label{LtildeBB}
\|\mathcal{L}_{\widetilde{B},B}\|\leq2\|B\|+2|a_2c_1|\|R\|^2\leq2\max(|b_1|,2|b_2|)+4|a_2c_1|\leq2+4=6,
\end{equation}
and the following lemma deduced from (\ref{LWAtendsVersId}).
\begin{Lemma}\label{LWABounded}
For $h$ small enough, it holds for all
$X\in\mathcal{M}_{(J+1)^2}(\mathbb{R})$  that
$$
\displaystyle\frac{1}{2}\|X\|\leq(1-C(\alpha)h^{2s})\|X\|\leq\|\mathcal{L}_{W,A}(X)\|\leq(1+C(\alpha)h^{2s})\|X\|\leq\displaystyle\frac{3}{2}\|X\|.
$$
\end{Lemma}
As a result, (\ref{LyapunovStability3}) yields that
\begin{equation}\label{LyapunovStability3-1}
\displaystyle\frac{1}{2}\|U^{n+1}\|\leq6\|U^{n}\|+2\|V^n\|+\displaystyle\frac{3}{2}\|U^{n-1}\|+2(\|F^{n-1}\|+\|F^n\|).
\end{equation}
For $n=0$, this implies that
\begin{equation}\label{LyapunovStability3-2}
\|U^{1}\|\leq12\|U^{0}\|+4\|V^0\|+3\|U^{-1}\|+4(\|F^{-1}\|+\|F^0\|).
\end{equation}
Using the discrete initial condition
$$
U^0=U^{-1}+l\varphi.
$$
Here we identify the function $\varphi$ to the matrix whom
coefficients are $\varphi_{j,m}=\varphi(x_j,y_m)$. We obtain
\begin{equation}\label{U-1Bounds}
\|U^{-1}\|\leq\|U^0\|+l\|\varphi\|.
\end{equation}
Observing that
$$
F^{-1}_{j,m}=F(U^{-1}_{j,m})=(U^{0}_{j,m}-l\varphi_{j,m})^2,
$$
it results that
$$
|F^{-1}_{j,m}|\leq|U^{0}_{j,m}|^2+2l|\varphi_{j,m})|.|U^{0}_{j,m}|+l^2|\varphi_{j,m}|^2
$$
and consequently,
\begin{equation}\label{F-1Bounds}
\|F^{-1}\|\leq\|U^{0}\|^2+2l\|\varphi\|.\|U^{0}\|+l^2\|\varphi\|^2.
\end{equation}
Hence, equation (\ref{LyapunovStability3-2}) yields that
\begin{equation}\label{LyapunovStability3-3}
\|U^{1}\|\leq(15+8l\|\varphi\|)\|U^{0}\|+4\|V^0\|+8\|F^0\|+3l\|\varphi\|+4l^2\|\varphi\|^2.
\end{equation}
Now, the Lyapunov criterion for stability states exactly that
\begin{equation}\label{LyapunovStability1}
\forall\,\,\varepsilon>0,\,\exists\,\eta>0\,\,\,s.t;\,\,\|(U^0,V^0)\|
\leq\eta\,\,\Rightarrow\,\,\|(U^n,V^n)\|\leq\varepsilon,\,\,\forall\,n\geq0.
\end{equation}
For $n=1$ and $\|(U^1,V^1)\|\leq\varepsilon$, we seek an $\eta>0$
for which $\|(U^0,V^0)\|\leq\eta$. Indeed, using
(\ref{LyapunovStability3-3}), this means that, it suffices to find
$\eta$ such that
\begin{equation}\label{LyapunovStability3-4}
8\eta^2+(19+8l\|\varphi\|)\eta+3l\|\varphi\|+4l^2\|\varphi\|^2-\varepsilon<0.
\end{equation}
Choosing $l$ small enough, we obtain a discriminant
$$
\Delta=(19+8l\|\varphi\|)^2-32(3l\|\varphi\|+4l^2\|\varphi\|^2-\varepsilon)>0
$$
and a zero
$\eta_1=\displaystyle\frac{\sqrt{\Delta}-(19+8l\|\varphi\|)}{16}>0$.
Consequently, it suffices to choose $\eta\in]0,\eta_1[$ to obtain
(\ref{LyapunovStability3-4}). Finally,
(\ref{LyapunovStability3-3}) yields that
$\|U^1\|\leq\varepsilon$.\\
Now, equation (\ref{LyapunovStability4}), for $n=0$, implies that
\begin{equation}\label{LyapunovStability4-1}
\|V^{1}\|\leq\,A(l,h,\varphi)\|U^{0}\|^2+B(l,h,\varphi)\|U^{0}\|+C(l,h,\varphi)+16|c_2|\|V^{0}\|,
\end{equation}
where
$$
A(l,h,\varphi)=3+32|c_2|,
$$
$$
B(l,h,\varphi)=4\left(|c_1|+8|c_2|(2+l\varphi)+l\varphi+\displaystyle\frac{1}{h^2}\right),
$$
and
$$
C(l,h,\varphi)=2(1+8|c_2|)l^2\|\varphi\|^2+4l(4|c_2|+\displaystyle\frac{1}{h^2})\|\varphi\|.
$$
Choosing $\|(U^0,V^0)\|\leq\eta$, it suffices to study the
inequality
\begin{equation}\label{V1Bound}
A(l,h,\varphi)\eta^2+\left(B(l,h,\varphi)+16|c_2|\right)\eta+C(l,h,\varphi)-\varepsilon\leq0.
\end{equation}
Its discriminant satisfies
$$
\Delta\sim\displaystyle\frac{16}{h^4}\left(1+20\alpha+|1-2\alpha|\right)^2
-\displaystyle\frac{128\alpha}{h^2}\left((4+16\alpha)\|\varphi\|h^s-\varepsilon\right)>0
$$
for $h$ small enough. The zero
$\eta_1'=\displaystyle\frac{\sqrt{\Delta}-(B(l,h,\varphi)+16|c_2|)}{A(l,h,\varphi)}>0$.
As a consequence, for $\eta\in]0,\eta_1'[$ we obtain
$\|V^1\|\leq\varepsilon$. Finally, for $\eta\in]0,\eta_0[$ with
$\eta_0=\min(\eta_1,\eta_1')$, we obtain
$\|(U^1,V^1)\|\leq\varepsilon$ whenever $\|(U^0,V^0)\|\leq\eta$.
Assume now that the $(U^k,V^k)$ is bounded for $k=1,2,\dots,n$ (by
$\varepsilon_1$) whenever $(U^0,V^0)$ is bounded by $\eta$ and let
$\varepsilon>0$. We shall prove that it is possible to choose
$\eta$ satisfying $\|(U^{n+1},V^{n+1})\|\leq\varepsilon$. Indeed,
from (\ref{LyapunovStability3-1}), we have
\begin{equation}\label{LyapunovStability3Ordren-1}
\|U^{n+1}\|\leq19\varepsilon_1+8\varepsilon_1^2.
\end{equation}
So, one seeks, $\varepsilon_1$ for which
$8\varepsilon_1^2+19\varepsilon_1-\varepsilon\leq0$. The
discriminant is $\Delta=361+32\varepsilon$, leading to one
positive zeros
$\varepsilon_1=\displaystyle\frac{\sqrt{361+32\varepsilon}-19}{16}$.
Then $\|U^{n+1}\|\leq\varepsilon$ whenever
$\|(U^k,V^k)\|\leq\varepsilon_1$, $k=1,2,\dots,n$. Next, using
(\ref{LyapunovStability4}) and (\ref{LyapunovStability3Ordren-1}),
we have
\begin{equation}\label{LyapunovStability3Ordren-2}
\|V^{n+1}\|\leq\left(4|c_1|+80|c_2|+1\right)\varepsilon_1+\left(32|c_2|+2\right)\varepsilon_1^2.
\end{equation}
So, it suffices as previously to choose $\varepsilon_1$ such that
$$
\left(32|c_2|+2\right)\varepsilon_1^2+\left(4|c_1|+80|c_2|+1\right)\varepsilon_1-\varepsilon\leq0.
$$
The discriminant is
$\Delta=(4|c_1|+80|c_2|+1)^2+4(32|c_2|+2)\varepsilon$, leading to
one positive zeros
$\varepsilon_1'=\displaystyle\frac{\sqrt{\Delta}-(4|c_1|+80|c_2|+1)}{2(32|c_2|+2)}$.
Then $\|U^{n+1}\|\leq\varepsilon$ whenever
$\|(U^k,V^k)\|\leq\varepsilon_1'$, $k=1,2,\dots,n$. Next, it holds
from the recurrence hypothesis for
$\varepsilon_0=\min(\varepsilon_1,\varepsilon_1')$, that there
exists $\eta>0$ for which $\|(U^0,V^0)\|\leq\eta$ implies that
$\|(U^k,V^k)\|\leq\varepsilon_0$, for $k=1,2,\dots,n$, which by
the next induces that $\|(U^{n+1},V^{n+1})\|\leq\varepsilon$.
\section{Numerical implementations}
The initial data becomes
\begin{equation}\label{CID1}
U_{j,m}^0=u(x_j,y_m,t_0)=u_0(x_j,y_m),\quad(j,m)\in\,I,
\end{equation}
\begin{equation}\label{CID2}
U_{j,m}^1=\displaystyle\frac{1}{2}\left(2u_0(x_j,y_m)+l^2\Delta\,u_0(x_j,y_m)+(v_0)_{xx}(x_j,y_m)-2lg(x_j,y_m)\right).
\end{equation}

\end{document}